\documentclass{article}
\usepackage{graphicx} % Required for inserting images
\usepackage[utf8]{inputenc}
\usepackage[T1]{fontenc}

\newcommand{\KP}{\mathsf{KP}}

\newcommand{\ZFC}{\mathsf{ZFC}}
\newcommand{\Ord}{\mathrm{Ord}}

\newcommand{\mytitle}{Recursive Analogues of Shrewdness and Subtlety, with Applications to Fine Structure}
\newcommand{\myauthor}{Jayde S. Massmann}

\font \rfont = cmr12 at 15pt
\font \authorfont = cmr12 at 12.5pt

\title{\rfont \mytitle}
\author{\authorfont \myauthor}
\date{\rfont \today}

\usepackage{hyperref}
\usepackage{comment}
\usepackage{proof}

\usepackage[utf8]{inputenc}
\usepackage{xcolor}
\usepackage{biblatex}
\usepackage{csquotes}
\usepackage{amsfonts}
\usepackage{amsmath}
\usepackage{cleveref}
\usepackage{amssymb}
\usepackage{amsthm}
\usepackage{parskip}
\usepackage[a4paper, total={6in, 8in}]{geometry}
\usepackage{titling}
\usepackage{tikz}
\usetikzlibrary{cd}
\usetikzlibrary{matrix,arrows}
\usetikzlibrary{decorations.pathreplacing,calc,arrows.meta}
\usepackage[normalem]{ulem}
\usepackage[english]{babel}

\usepackage{blindtext}

\usepackage{bbm}

\theoremstyle{definition}
\newtheorem{definition}{Definition}[section]

\theoremstyle{plain}
\newtheorem{theorem}[definition]{Theorem}

\theoremstyle{plain}

\theoremstyle{plain}
\newtheorem{proposition}[definition]{Proposition}

\theoremstyle{plain}
\newtheorem{corollary}[definition]{Corollary}

\theoremstyle{remark}

\theoremstyle{remark}

\theoremstyle{remark}

\theoremstyle{plain}
\newtheorem{lemma}[definition]{Lemma}

\theoremstyle{plain}

\theoremstyle{plain}

\begin{document}

\maketitle

\begin{abstract}
We share both recent and older, well-known results regarding the notions of stable ordinals and shrewd cardinals. We then argue that $\Sigma_2$-nonprojectible ordinals may be considered as recursive analogues to subtle cardinals, a highly combinatorial type of cardinal related to Jensen's fine structure, due to the latter possessing a characterisation in terms of shrewdnesss.
\end{abstract}

\section{Introduction}

The notion of a shrewd cardinal was introduced by Rathjen in \cite{rathjen}. Shrewd cardinals offer an alternative transfinite extension of indescribability to $\xi$-indescribability (for $\xi \geq \omega$), since the latter suffers from the following two unfavourable properties:

\begin{enumerate}
    \item If $\kappa$ is $\xi$-indescribable and $\xi' < \xi$, $\kappa$ need not be $\xi'$-indescribable.
    \item $\kappa$ cannot be $\kappa$-indescribable.
\end{enumerate}

Which imply that the hierarchy is not so nice and stops abruptly, being unable to reach, say, ``hyper-indescribability''. Shrewdness suffers from neither of these flaws, and therefore can be considered as a linear hierarchy of large cardinal axioms bridging the gap between weak compactness or finite stages of indescribability, and more powerful combinatorial or elementary embedding-based large cardinals in the literature. Let us give their definition:

\begin{definition}
\label{ShrewdCardinal}
Let $\eta > 0$. A cardinal $\kappa$ is $\eta$-shrewd iff, for all $P \subseteq V_\kappa$ and every formula $\varphi(x)$, if $(V_{\kappa+\eta}, \in, P) \models \varphi(\kappa)$, then there exist $0 < \kappa_0, \eta_0 < \kappa$ so that $(V_{\kappa_0+\eta_0}, \in, P \cap V_{\kappa_0}) \models \varphi(\kappa_0)$.
\end{definition}

By slightly modifying the definition, one can obtain various natural strengthenings of these cardinals, as well as delineate the consistency strength hierarchy of shrewdness more finely. In particular, we will consider the definition of two additional parameters: an arbitrary class $\mathcal{A}$ and a set $\mathcal{F}$ of formulae.

The definition of a $\mathcal{A}$-$\eta$-$\mathcal{F}$-shrewd cardinal is obtained from the definition of an $\eta$-shrewd cardinal by:

\begin{enumerate}
    \item Restricting the formula $\varphi$ to be an element of $\mathcal{F}$.
    \item Adding an additional predicate for $\mathcal{A}$ in both the hypothesis and the conclusion. Namely, $(V_{\kappa+\eta}, \in, P) \models \varphi(\kappa)$ becomes $(V_{\kappa+\eta}, \in, P, \mathcal{A} \cap V_{\kappa+\eta}) \models \varphi(\kappa)$, and similarly $(V_{\kappa_0+\eta_0}, \in, P \cap V_{\kappa_0}) \models \varphi(\kappa_0)$ becomes $(V_{\kappa_0+\eta_0}, \in, P \cap V_{\kappa_0}, \mathcal{A} \cap V_{\kappa_0+\eta_0}) \models \varphi(\kappa_0)$.
\end{enumerate}

The formula therefore may now contain two additional predicate symbols, one for first-order variables (representing $P$) and one for the highest-order variables (representing $\mathcal{A}$). Therefore, it makes sense to consider natural $\mathcal{F}$'s which are closed under addition of predicates. The most natural $\mathcal{F}$'s to consider are therefore the classes $\Pi_n$ and $\Sigma_n$ of the Lévy hierarchy.

Note that if $\mathcal{A} = V$, one may eliminate the second part and so consider simply $\eta$-$\mathcal{F}$-shrewdness. Dually, if $\mathcal{F}$ consists of all formulae, one may eliminate the first part and so consider simply $\mathcal{A}$-$\eta$-shrewdness. By eliminating both, one arrives where we started -- $\eta$-shrewdness.

It is now an immediate observation that $\kappa$ is $\eta$-shrewd iff it is $\eta$-$\Pi_n$-shrewd for all $n$. This is since any formula is $\Pi_n$ for some $n$, by writing in prenex normal form. This allows us to connect back finite stages of shrewdness to finite stages of indescribability:

\begin{proposition}
\label{ShrewdToIndesc}
(Folklore) Assume $0 \leq n, m < \omega$ and $n \neq 0$. Then $\kappa$ is $n$-$\Pi_m$-shrewd iff $\kappa$ is $\Pi^n_m$-indescribable. In particular, $\kappa$ is $n$-shrewd iff $\kappa$ is $\Pi^{n+1}_0$-indescribable.
\end{proposition}

The proof that, if $\kappa$ is $\eta$-shrewd and $\delta < \eta$, then $\kappa$ is $\delta$-shrewd (which implies that, say, if $\kappa$ is $\eta$-shrewd and $\eta \geq \omega$, then $\kappa$ is totally indescribable) can be found in \cite{rathjen2}. 

Subtlety is a combinatorial principle introduced by Jensen and Kunen relating to the former's analysis of the combinatorial and fine structural properties of the constructible universe $L$. While originally formulated in a manner similarly to the $\diamond$ principle, it also has a characterisation in terms of shrewdness. This characterisation is quite reminiscent of Woodin and Vopěnka cardinals. Let us first give the original definition.

\begin{definition}
\label{SubtleDef}
A cardinal $\kappa$ is called subtle iff, for any sequence $\langle S_\alpha: \alpha < \kappa \rangle$ satisfying $S_\alpha \subseteq \alpha$ for all $\alpha < \kappa$, and club $C$ in $\kappa$, there are $\beta, \delta \in C$ so that $\beta < \delta$ and $S_\delta \cap \beta = S_\beta$.
\end{definition}

In essence, any sequence $\langle S_\alpha: \alpha < \kappa \rangle$ satisfying $S_\alpha \subseteq \alpha$ for all $\alpha < \kappa$ can be made to ``cohere'' club many times. Some properties of subtle cardinals can be found in \cite{jensen}, e.g. if $\kappa$ is subtle then $\diamond_\kappa$ holds. For the following characterisation, say $\kappa$ is $\mathcal{A}$-$< \pi$-shrewd iff it is $\mathcal{A}$-$\eta$-shrewd for all $\eta < \pi$.

\begin{theorem}
\label{SubtleFromShrewd}
Let $\pi$ be an inaccessible cardinal. Then TFAE:

\begin{enumerate}
    \item $\pi$ is subtle.
    \item For any $\mathcal{A} \subseteq V_\pi$, the set of $\kappa$ which are $\mathcal{A}$-$< \pi$-shrewd is stationary in $\pi$.
\end{enumerate}
\end{theorem}

\begin{proof}
The forward direction can be found in \cite{rathjen2}. For the sake of completeness, we will give it here. So, assume $\pi$ is subtle. It is known that $\pi$ is strongly inaccessible, and thus $|V_\pi| = \pi$. Let $F: V_\pi \to \pi$ be a bijection. And since strongly inaccessible cardinals are limits of cardinals $\kappa < \pi$ so that $|V_\kappa| = \kappa$, one may pick $F$ in a way so that the set of $\kappa < \pi$ with $F''V_\kappa = \kappa$ is club. This club is denoted $C_F$. Now pick an arbitrary club $B$ in $\pi$. We aim to show that $B$ contains a $\mathcal{A}$-$< \pi$-shrewd cardinal. Since $\pi$ is uncountable regular, the intersection of two clubs in $\pi$ is also club, so we may assume without loss of generality that $B \subseteq C_F$, replacing $B$ with $B \cap C_F$ if this isn't the case. In particular, all elements of $B$ can be taken to be cardinals.

Now for a contradiction assume there is no cardinal $\kappa \in B$ so that $\kappa$ is $\mathcal{A}$-$< \pi$-shrewd. One can use the axiom of choice to pick, for all $\kappa \in B$, a $\sigma_\kappa$ witnessing that $\kappa$ is not $\mathcal{A}$-$< \pi$-shrewd (i.e. for each $\kappa \in B$, $\kappa$ is not $\mathcal{A}$-$\sigma_\kappa$-shrewd). We may once again assume, without loss of generality, that $\kappa < \sigma_\kappa$ for all $\kappa \in B$, by the fact that if $\kappa$ is $\mathcal{A}$-$\eta$-shrewd and $\delta < \eta$, then $\kappa$ is $\mathcal{A}$-$\delta$-shrewd (and so, by contrapositive, if $\kappa$ is not $\mathcal{A}$-$\eta$-shrewd for some $\eta \leq \kappa$, then $\kappa$ is not $\mathcal{A}$-$\eta$-shrewd for any $\eta > \kappa$). We can extend $\sigma$ to a total map with domain $\kappa$ by simply putting $\sigma_\rho = \rho$ for $\rho \notin B$.

Let $E$ be the set of $\rho \in B$ so that $\sigma_\nu < \rho$ whenever $\nu < \rho$. For example, $\min(B) \in E$, since $\nu < \min(B)$ implies $\nu \notin B$ and whence $\sigma_\nu = \nu$. A simple argument shows that $E$ is club in $\pi$ as well. It's easy to see that, if $\kappa_0, \kappa_1 \in E$ and $\kappa_0 < \kappa_1$, then $\sigma_{\kappa_0} < \kappa_1$ and so $\kappa_0$ is not $\mathcal{A}$-$\kappa_1$-shrewd. Therefore, for $\kappa \in E$, let $\kappa^*$ be the minimal element of $E$ above $\kappa$. Note that, by our hypothesis that all elements of $B$ are cardinals, we have $\kappa + \kappa^* = \kappa$. Now since $\kappa$ is never $\mathcal{A}$-$\kappa^*$-shrewd, we may let $\varphi_\kappa$ and $P_\kappa$ witness this, i.e. be so that $(V_{\kappa^*}, \in, P_\kappa, \mathcal{A} \cap V_{\kappa^*}) \models \varphi_\kappa(\kappa)$, but whenever $0 < \nu, \delta < \kappa$, we have $(V_{\nu+\delta}, \in, P_\kappa \cap V_\nu, \mathcal{A} \cap V_{\nu+\delta}) \models \neg \varphi_\kappa(\nu)$.

Now let $\theta_\kappa(u, v)$ be the formula $v \in \Ord \land \exists \xi > v (u \subseteq V_\xi \land (V_\xi, \in, u, \mathbf{P} \cap V_\xi) \models \varphi_\kappa(v))$. Here, $\mathbf{P}$ is the predicate (i.e. definable class) which is interpreted as $\mathcal{A}$. Now, if $\kappa^* < \rho$, then $(V_\rho, \in, \mathcal{A} \cap V_\rho) \models \theta_\kappa(P_\kappa, \kappa)$, witnessed by $\xi = \kappa^*$. Let $\vartheta_\kappa(v)$ be the formula in the language with both additional predicate symbols which states that $\theta_\kappa(u, v)$ holds where $u$ is the subclass of the first-order domain of all objects satisfying the first predicate. Therefore, $\kappa^* < \rho$ implies $(V_\rho, \in, P_\kappa, \mathcal{A} \cap V_\rho) \models \vartheta_\kappa(\kappa)$. Furthermore, for all $0 < \mu < \kappa$, we have $(V_\kappa, \in, P_\kappa \cap V_\mu, \mathcal{A} \cap V_\kappa) \models \neg \vartheta_\kappa(\mu)$ by our assumption that whenever $0 < \nu, \delta < \kappa$, we have $(V_{\nu+\delta}, \in, P_\kappa \cap V_\nu, \mathcal{A} \cap V_{\nu+\delta}) \models \neg \varphi_\kappa(\nu)$.

Now let $E_\infty$ be the set of ordinals below $\pi$ which are limit points of $E$. The upshot is that, whenever $\kappa, \rho \in E_\infty$ and $\kappa < \rho$, then $(V_\rho, \in, P_\kappa, \mathcal{A} \cap V_\rho) \models \vartheta_\kappa(\kappa)$ and, whenever $0 < \mu < \kappa$, $(V_\kappa, \in, P_\kappa \cap V_\mu, \mathcal{A} \cap V_\kappa) \models \neg \vartheta_\kappa(\mu)$. Now we use subtlety to derive a contradiction. Let $\langle \psi_n: n < \omega \rangle$ enumerate all formulae with one free variable in our language with two predicate symbols. Also, for $\kappa \in E_\infty$, let $\kappa^{**}$ denote the minimal element of $E_\infty$ above $\kappa$. Then define $P_\kappa^*$ as a particular coding of $P_\kappa$, namely $P_\kappa^* = (F''P_\kappa \cap (\kappa \setminus \omega)) \cup \{3n: n \in F''P_\kappa \cap \omega\} \cup \{3n+1: (V_{\kappa^{**}}, \in, P_\kappa, \mathcal{A} \cap V_{\kappa^{**}}) \models \psi_n(\kappa)\} \cup \{3n+2: (V_{\kappa^{**}}, \in, P_\kappa, \mathcal{A} \cap V_{\kappa^{**}}) \models \neg \psi_n(\kappa)\}$. For $\alpha \notin E_\infty$, set $P_\alpha^* = \alpha$. 

Since $E_\infty$ is the set of limit points of a particular club and is therefore seen to also be a club, and it is also immediate that $P_\alpha^* \subseteq \alpha$ for all $\alpha < \delta$, we can use subtlety of $\pi$ to find two $\beta, \gamma \in E_\infty$ so that $\beta < \gamma$ and $P_\gamma^* \cap \beta = P_\beta^*$. By the way we defined $P_\kappa^*$ for $\kappa \in E_\infty$, the fact that $F$ is injective, and that $E_\infty \subseteq C_F$, it follows that $P_\gamma \cap V_\beta = P_\beta$. By our previous remark that $\kappa, \rho \in E_\infty$ and $\kappa < \rho$ implies $(V_\rho, \in, P_\kappa, \mathcal{A} \cap V_\rho) \models \vartheta_\kappa(\kappa)$, it follows that $(V_{\gamma^{**}}, \in, P_\gamma, \mathcal{A} \cap V_{\gamma^{**}}) \models \vartheta_\gamma(\gamma)$. Combining this with $P_\gamma^* \cap \beta = P_\beta^*$ and the definition of $P_\kappa^*$ for $\kappa \in E_\infty$, we obtain $(V_{\beta^{**}}, \in, P_\gamma \cap V_\beta, \mathcal{A} \cap V_{\beta^{**}}) \models \vartheta_\gamma(\beta)$. Another application of the precise way $P_\kappa^*$ is defined for $\kappa \in E_\infty$, together with $P_\gamma^* \cap \beta = P_\beta^*$, we get $(V_\gamma, \in, P_\gamma \cap V_\beta, \mathcal{A} \cap V_\gamma) \models \vartheta_\gamma(\beta)$. But this contradicts our assumption that, whenever $0 < \mu < \kappa$, $(V_\kappa, \in, P_\kappa \cap V_\mu, \mathcal{A} \cap V_\kappa) \models \neg \vartheta_\kappa(\mu)$.

Now, for the backwards direction, we utilize Rathjen's notion of reducibility. Assume that $\pi$ is inaccessible and, for any $\mathcal{A} \subseteq V_\pi$, the set of $\kappa$ which are $\mathcal{A}$-$< \pi$-shrewd is stationary in $\pi$. Via Appendix 7.8 of \cite{rathjen}, we see that any $\mathcal{A}$-$< \pi$-shrewd cardinal below $\pi$ is $\mathcal{A}$-$< \pi$-reducible. Now let $C$ be a club in $\pi$ and $\vec{S} = \langle S_\alpha: \alpha < \pi \rangle$ satisfy $S_\alpha \subseteq \alpha$ for all $\alpha < \pi$. It actually suffices for $\vec{S}$ to be partial with domain $C$, i.e. $S_\alpha$ is only defined for $\alpha \in C$. This will be relevant later, and does not affect the large cardinal axiom as we only consider $S_\beta$ for $\beta \in C$. Note then that $\vec{S} \subseteq V_\pi$. By stationarity, there is $\kappa \in C$ which is $\vec{S}$-$< \pi$-reducible. We will be interested in the specific case that $\kappa$ is $\vec{S}$-$\kappa + 2$-reducible. Now there is some $0 < \kappa_0 < \eta_0 < \kappa + 2$ so that $\langle V_{\eta_0}, \in, \kappa_0, \vec{S}|\eta_0 \rangle$ is elementarily equivalent to $\langle V_{\kappa+2}, \in, \kappa, \vec{S}|(\kappa + 2) \rangle$. We first claim that we must have $\eta_0 = \kappa_0+2$: let $\varphi$ be a formula formalizing ``$c+1$ exists and is the largest ordinal'', where $c$ is a constant symbol representing $\kappa_0$ or $\kappa$. Since $\langle V_{\kappa+2}, \in, \kappa \rangle \models \varphi$, we have $\langle V_{\eta_0}, \in, \kappa_0 \rangle \models \varphi$, and so $\kappa_0 + 1$ is the largest ordinal of $V_{\eta_0}$. Thus $\eta_0 = \kappa_0 + 2$.

Utilizing elementary equivalence again, but this time exploiting the predicate, we obtain, for any $\alpha < \kappa_0$, $\langle V_{\eta_0}, \in, \kappa_0, \vec{S}|(\kappa_0 + 2) \rangle \models \alpha \in S_{\kappa_0}$ iff $\langle V_{\kappa+2}, \in, \kappa, \vec{S}|(\kappa + 2) \rangle \models \alpha \in S_\kappa$. Thus, $S_\kappa \cap \kappa_0 = S_{\kappa_0}$, as desired. It remains to show that $\kappa, \kappa_0 \in C$. We already have $\kappa \in C$ by hypothesis, and

\begin{equation}
\begin{split}
\kappa_0 \in C & \iff \kappa_0 \in \operatorname{dom}(\vec{S}) \\ & \iff \kappa_0 \in \operatorname{dom}(\vec{S}|(\kappa_0 + 2)) \\ & \iff \langle V_{\kappa_0 + 2}, \in, \kappa_0, \vec{S}|(\kappa_0 + 2) \rangle \models ``\kappa_0 \in \operatorname{dom}(\vec{S}|(\kappa_0 + 2)) \\ & \iff \langle V_{\kappa + 2}, \in, \kappa, \vec{S}|(\kappa_0 + 2) \rangle \models ``\kappa \in \operatorname{dom}(\vec{S}|(\kappa + 2)) \\ & \iff \kappa \in C
\end{split}
\end{equation}

which is true.
\end{proof}

In particular, if $\pi$ is subtle, then for any $\eta < \pi$ there is an $\kappa < \pi$ which is $\eta$-shrewd. Of course, the full characterisation we gave is stronger. Note that $\pi$ itself need not be $1$-shrewd: this is since any $1$-shrewd cardinal is $\Pi^2_0$-indescribable, but subtlety is describable by a $\Pi^1_1$ definition and whence the least subtle cardinal isn't even weakly compact.

\section{Recursive analogues and stable ordinals}

In this brief note, we shall show that $\Sigma_2$-nonprojectible ordinals may be considered as so-called ``recursive analogues'' of subtle cardinals. These are strengthenings of the notion of a nonprojectible ordinal, which was a direct byproduct of Jensen's analysis of the fine structure of $L$, thereby also being perhaps genealogically related to subtlety as well. This notion was first shared in print in the famous article \cite{jensen2}, under the name of a strongly admissible ordinal. The term ``nonprojectibility'' arises from the fact that $\alpha$ is nonprojectible iff, whenever $\eta < \alpha$, $X \subseteq \eta$ and $f: X \to \alpha$ is surjective, then $f$ is not $\Sigma_1$-definable in $L_\alpha$ with parameters. Equivalently, the ordinals $\eta < \alpha$ so that $L_\eta \prec_{\Sigma_1} L_\alpha$ are unbounded, which is the definition we shall use. The original definition was as ordinals $\alpha$ so that $L_\alpha$ satisfies Kripke-Platek set theory $\KP$ with $\Sigma_1$-separation, which also is equivalent to ordinals $\alpha$ so that $L_\alpha \cap \mathcal{P}(\omega)$ is a model of $\Pi^1_2$-comprehension, where the latter additionally required $\alpha < \omega_1^L$. See \cite{barwise} and \cite{simpson} for proofs of these characterisations. $\Sigma_2$-nonprojectibility can be obtained from any of these characterisations by going ``one level up'' in complexity: e.g. not having $\Sigma_2$-definable surjective mappings, having $\Sigma_2$-elementary substructures, satisfying $\Sigma_2$-separation or modelling $\Pi^1_3$-comprehension. We will once again primarily be dealing with the characterisation via elementary substructures.

The notion of a recursive analogue is relevant to $\alpha$-recursion theory, in which one generalizes theorems in classical recursion and computational complexity theory (where primitive objects are hereditarily finite sets, e.g. natural numbers and finite binary sequences) to theorems about $L_\alpha$ where $\alpha$ is an admissible ordinal. For example, for $\alpha > \omega$, hyperarithmetic reals may be treated as primitive objects as they arise in $L_{\omega_1^{\mathrm{CK}}} \cap \mathcal{P}(\omega)$. In essence, a recursive analogue of a large cardinal is a notion describing large countable ordinals, which may behave in a way that mimics the way that large cardinals act. This also makes them useful in proof theory, as the existence of recursive analogues of large cardinals is actually provable in $\ZFC$, unlike the existence of large cardinals, and so one may define ordinal representation systems in $\ZFC$ without requiring large cardinal hypotheses. For example, an ordinal analysis of $\KP + ``$ every set is contained in a standard transitive model of $\KP$'', also denoted $\mathsf{KPI}$, has traditionally required the additional assumption of the existence of a weakly inaccessible cardinal, but this could be eliminated by replacing weakly inaccessible cardinals with their recursive analogue -- recursively inaccessible ordinals.

Recursively regular ordinals are typically considered to be admissible ordinals (and whence recursively inaccessible ordinals are precisely the admissible limits of admissible ordinals), because of the similarity between the following two characterisations:

\begin{lemma}
A cardinal $\kappa$ is regular iff, for any function $f: \kappa \to \kappa$, there is an $\alpha < \kappa$ so that $f''\alpha \subseteq \alpha$. An ordinal $\kappa$ is admissible iff the above holds when $f$ is restricted to be $\Delta_1$-definable in $L_\kappa$ with parameters.
\end{lemma}

Unfortunately, this substitution means that slightly more heavy lifting is required, as one needs to verify that the desired projection functions are sufficiently definable. The notions of $\Delta_1$- and $\Sigma_1$-definability in $L_\alpha$ with parameters arises often, so we adopt the following abbreviative convention from generalized recursion theory: $A \subseteq L_\alpha$ is $\alpha$-recursive (resp. $\alpha$-recursively enumerable) iff it is $\Delta_1$-definable (resp. $\Sigma_1$-definable) in $L_\alpha$ with parameters, i.e. $\Delta_1(L_\alpha)$ (resp. $\Sigma_1(L_\alpha)$. Note that any $\alpha$-recursively enumerable map $f: L_\alpha \to L_\alpha$ is already $\alpha$-recursive, since $f(x) \neq y$ is equivalent to $\exists z (f(x) = z \land z \neq y)$ and $\Sigma_1$ formulae are closed under conjunction and existential quantification. Back to the topic at hand, let us define the notion of a $\xi$-$\Pi_n$-reflecting ordinal.

\begin{definition}
\label{ReflectingOrdinal}
Let $A$ be a class of ordinals. For $\xi > 0$, an ordinal $\alpha$ is called $\xi$-$\Pi_n$-reflecting onto $A$ iff, for every $\Pi_n$-formula $\varphi(x,\vec{y})$, and for all $\vec{b} \in L_\alpha$, if $L_{\alpha+\xi} \models \varphi(\alpha, \vec{b})$, then there exist $\alpha_0, \xi_0 < \alpha$ so that $\xi_0 > 0$, $\alpha_0 \in A$, $\vec{b} \in L_{\alpha_0}$ and $L_{\alpha_0+\xi_0} \models \varphi(\alpha_0, \vec{b})$.

$\alpha$ is called $\xi$-$\Pi_n$-reflecting iff it is $\xi$-$\Pi_n$-reflecting onto $\Ord$. And $\alpha$ is $\Pi_n$-reflecting iff it is ``$0$-$\Pi_n$-reflecting'', which means that the parameter for $\alpha$ is eliminated (since $\alpha \notin L_\alpha)$), i.e. for  all $\vec{b} \in L_\alpha$, if $L_\alpha \models \varphi(\vec{b})$, then there exists $\alpha_0 < \alpha$ so that $\alpha_0 \in A$, $\vec{b} \in L_{\alpha_0}$ and $L_{\alpha_0} \models \varphi(\vec{b})$.
\end{definition}

It is known that $\alpha$ is $\Pi_2$-reflecting iff $\alpha > \omega$ and it is admissible. See \cite{arai} for a proof of this. Therefore, it may be argued that $\Pi_2$-reflecting ordinals serve as a countable analogue of uncountable regular ordinals, as already mentioned. It is known that $\Pi_3$-reflecting ordinals are $\Pi_2$-reflecting onto the class of $\Pi_2$-reflecting ordinals, and much more. Hence, it was argued by Richter and Aczel in \cite{richter} that $\Pi_{n+2}$-reflecting ordinals should serve as recursive analogues to $\Pi^1_n$-describable cardinals for $n > 0$. In general, for $\xi > 0$, $\xi$-$\Pi_n$-reflecting ordinals should serve as recursive analogues to $\xi+1$-$\Pi_n$-shrewd ordinals. A large focus has become stability, which links back to our previous mention of nonprojectibility:

\begin{definition}
\label{StableOrdinal}
Say $\alpha$ is $\xi$-stable iff $L_\alpha \prec_{\Sigma_1} L_\xi$, where $\prec_{\Sigma_1}$ denotes the relation of being a $\Sigma_1$-elementary substructure.
\end{definition}

It is easy to see that if $\alpha$ is $\alpha+1$-stable, then it is $\Pi_n$-reflecting for all $n < \omega$. Namely, let $\varphi(x)$ be an arbitrary formula, and $b \in L_\alpha$. Assume $L_\alpha \models \varphi(b)$. Then $L_{\alpha+1}$ satisfies ``there is an ordinal $\beta$ so that $b \in L_\beta$ and $L_\beta \models \varphi(b)$'', which can be written in $\Sigma_1$ form. Thus, by $L_\alpha \prec_{\Sigma_1} L_{\alpha+1}$, $L_\alpha$ satisfies the same thing, and so the $\beta$ witnessing this satisfies $b \in L_\beta$ and $L_\beta \models \varphi(b)$. It turns out that the converse is also true, and a more general result holds for $\xi$-$\Pi_n$-reflection.

\begin{lemma}
\label{StableVersusReflecting}
(Folklore) For an ordinal $\alpha$, $\alpha$ is $\alpha+\xi$-stable iff it is $\gamma$-$\Pi_n$-reflecting for all $n$ and $\gamma < \xi$.
\end{lemma}

Using Theorem \ref{SubtleFromShrewd} and Lemma \ref{StableVersusReflecting}, this motivates the definition of a recursively subtle ordinal.

\begin{definition}
Let $\mathcal{A}$ be an arbitrary class. Say that $\alpha$ is $\mathcal{A}$-$\xi$-stable iff $\langle L_\alpha, \in, \mathcal{A} \cap L_\alpha \rangle \prec_{\Sigma_1} \langle L_\xi, \in, \mathcal{A} \cap L_\xi \rangle$. Now say an ordinal $\rho$ is recursively subtle iff, for any $\rho$-recursively enumerable $\mathcal{A} \subseteq L_\rho$, $\rho$ is $\Pi_2$-reflecting onto the set of $\kappa < \rho$ which are $\mathcal{A}$-$\rho$-stable.
\end{definition}

Note that it is a known result that if $\kappa$ is $\mathcal{A}$-$< \rho$-stable, it is already $\mathcal{A}$-$\rho$-stable, and so we did not deviate from the ``spirit'' of Theorem \ref{SubtleFromShrewd} too much. Actually, let us briefly state some well-known properties of stability:

\begin{lemma}
\label{StabilityFacts}
\begin{enumerate}
    \item If $\alpha < \beta < \gamma$ and $\alpha$ is $\gamma$-stable, then $\alpha$ is $\beta$-stable.
    \item If $\alpha$ is $\beta$-stable and $\beta$ is $\gamma$-stable, then $\alpha$ is $\gamma$-stable.
    \item If $\alpha$ is $\beta$-stable for all $\alpha \in A$, then either $\sup A = \beta$ or $\sup A$ is $\beta$-stable.
    \item Dually, if $\alpha$ is $\beta$-stable for all $\beta \in A$, then $\alpha$ is $\sup B$-stable.
    \item If $\alpha$ is $\alpha+2$-stable, then it is $\Pi_n$-reflecting on the class of ordinals $\xi$ which are $\xi+1$-stable, for all $n < \omega$.
    \item The least ordinal $\alpha$ that is a limit of ordinals $\xi$ which are $\xi+1$-stable is not itself even admissible.
\end{enumerate}
\end{lemma}

Proofs of the sublemmata here are either folklore or due to the anonymous Discord user C7X with whom I've had many insightful conversations.

\begin{proof}
(1) This follows from an easy upwards absoluteness argument.

(2) This is trivial.

(3) Assume that, for all $\alpha \in A$, we have $L_\alpha \prec_{\Sigma_1} L_\beta$. We aim to show that $L_{\sup A} \prec_{\Sigma_1} L_\beta$, assuming $\sup A < \beta$. In the case when $\sup A = \max A$, it is trivial, so we may assume, without loss of generality, that $\sup A$ is strictly greater than all elements of $A$. Let $x \in L_{\sup A}$ and $\varphi$ be an arbitrary $\Sigma_1$ formula. We aim to show that $L_{\sup A} \models \varphi(x)$ iff $L_\beta \models \varphi(x)$. The forwards direction follows by upwards absoluteness. For the converse direction, assume $L_\beta \models \varphi(x)$. Since $\sup A$ is a limit ordinal and $A$ is cofinal in $\sup A$, there is some $\alpha \in A$ so that $x \in L_\alpha$. By $L_\alpha \prec_{\Sigma_1} L_\beta$, we have $L_\alpha \models \varphi(x)$. Then $L_{\sup A} \models \varphi(x)$, once again by upwards absoluteness.

(4) This follows by a similar argument to (3).

(5) Let $\varphi$ be a $\Pi_n$ formula so that $L_\alpha \models \varphi$, and $x \in L_\alpha$. Then $L_{\alpha+2} \models \exists \xi (L_\xi \prec_{\Sigma_1} L_{\xi+1} \land x \in L_\xi \land L_\xi \models \varphi(x))$, with witness $\xi = \alpha$, which can be rendered in $\Sigma_1$ form. Thus $L_\alpha \models \exists \xi (L_\xi \prec_{\Sigma_1} L_{\xi+1} \land x \in L_\xi \land L_\xi \models \varphi(x))$, and so there is some $\xi < \alpha$ so that $\xi$ is $\xi+1$-stable, $x \in L_\xi$ and $L_\xi \models \varphi$. Thus $\xi \in A \cap \alpha$ where $A$ is the class of $\xi$ which are $\xi+1$-stable, and so $\xi$ witnesses $\Pi_n$-reflection onto $A$ in this instance.

(6) This proof uses Richter and Aczel's notion of $\Sigma_1$-collection; note that this is different to many other notions of $\mathcal{F}$-collection. Now for contradiction, assume that $\alpha$ is admissible. Let $\varphi(n, \gamma)$ be a formula that asserts that $\gamma$ is the $n$'th, starting from $n = 0$, ordinal so that $\gamma$ is $\gamma+1$-stable. This formula can be defined uniformly in $n$. Then apply $\Sigma_1$-collection to the formula $L_\alpha \models \forall n < \omega \exists \nu \varphi(n, \nu)$ to obtain ``for some $b \in L_\alpha$, we have $L_\alpha \models \forall n < \omega \exists \nu (b \models \varphi(n, \nu))$'', i.e. some set in $L_\alpha$ contains $\omega$ many ordinals $\gamma$ so that $\gamma$ is $\gamma+1$-stable. But this contradicts minimality of $\alpha$.
\end{proof}

Many more advanced results about stability, e.g. that it is a proper hierarchy (i.e. that for all $\gamma$, there is an $\alpha$ that is $\alpha+\gamma$-stable but not $\alpha+\gamma+1$-stable), are proved via reflection, e.g. showing that an ordinal $\alpha$ which is $\alpha+\gamma+1$-stable is $\Pi_n$-reflecting on the set of ordinals $\xi < \alpha$ which are $\xi+\gamma$-stable, for all $n$. Much of the work in \cite{levy} is helpful for formalizing reflection schemata in first-order set theory.

The reason why we used $\Pi_2$-reflection in our formulation of recursive subtlety is because, recalling the statement that $\Pi_{n+2}$-reflection serves as a recursive analogue of $\Pi^1_n$-indescribability, it is known that $\kappa$ is $\Pi^1_0$-indescribable onto $A$ iff it is strongly inaccessible and $A$ is stationary in $\kappa$, thus if $\rho$ is recursively inaccessible (which we shall show all recursively subtle ordinals are) it yields another adequate generalizaton. For example, an ordinal is considered to be recursively Mahlo iff it is $\Pi_2$-reflecting onto the set of admissible ordinals below, analogously to how a cardinal is Mahlo iff it is $\Pi^1_0$-indescribable onto the set of regular cardinals below.

As mentioned, our main theorem is that $\rho$ is recursively subtle iff it is $\Sigma_2$-nonprojectible. First, we shall state some results regarding stability and nonprojectibility, including some direct implications and size comparisons. Recall that $\rho$'s nonprojectibility is a $\Pi_3$ property of $L_\rho$. Therefore, there is already some behaviour-wise similarity between nonprojectibility and subtlety: the least subtle cardinal is much greater than the least $1$-shrewd cardinal, and is itself Mahlo but not weakly compact. Similarly, the least nonprojectible ordinal is much greater than the least $\alpha$ which is $\alpha+1$-stable, and is (as we shall prove in a moment) itself recursively Mahlo but not $\Pi_3$-reflecting. $\Sigma_2$-nonprojectibility has many of the same ``behaviours''.

That nonprojectible ordinals are recursively Mahlo follows from the following theorem and the fact that if $\eta$ is $\rho$-stable, it is very obviously admissible:

\begin{theorem}
\label{NonprojectiblePi2Ref}
If $\rho$ is nonprojectible, then $\rho$ is $\Pi_2$-reflecting onto the set of $\kappa < \rho$ so that $\kappa$ is $\rho$-stable.
\end{theorem}

This theorem came as a surprise, since previously $\Pi_2$-reflection is generally considered stronger than even iterated limit point taking. However, there are already other contexts in which $\Pi_2$-reflection may fail to imply iterated limit points (e.g. $\aleph_\omega$ is $\Pi_2$-reflecting onto the set of cardinals, despite not being a limit of limit cardinals).

\begin{proof}
Assume that $\rho$ is nonprojectible, $L_\rho \models \forall x \exists y \varphi(x,y,b)$, where $\varphi(x,y,p)$ is a $\Delta_0$-formula, and $b \in L_\rho$. Let $\beta$ be so that $\beta$ is $\rho$-stable and $b \in L_\beta$. Now let $x \in L_\beta$. Since $\beta$ is $\rho$-stable, there is some $y \in L_\beta$ so that $\varphi^{L_\beta}(x,y,b)$ since $L_\rho \models \exists y \varphi(x,y,b)$ and $L_\beta \prec_{\Sigma_1} L_\rho$. Now let $y \in L_\beta$ be so that $L_\beta \models \varphi(x,y,b)$. Since $\varphi$ is $\Delta_0$ and such formulae are absolute for transitive sets, it follows that $\varphi^{L_\rho}(x,y,b)$ iff $\varphi^{L_\beta}(x,y,b)$, thus $L_\beta \models \forall x \exists y \varphi(x,y,b)$. The desired result follows.
\end{proof}

Also let us briefly state some results regarding nonprojectibility.

\begin{theorem}
Say $\alpha$ is $\omega$-fold stable iff there is a map $f: \omega \to \Ord$ so that $f(0) = \alpha$ and, for all $i < \omega$, $f(i)$ is $f(i+1)$-stable.

\begin{enumerate}
    \item Assume $\alpha$ is $\omega$-fold stable, witnessed by $f$. Then $\sup\{f(i): i < \omega\}$ is nonprojectible.
    \item As a sort of converse, assume $\rho$ is nonprojectible, $\operatorname{cof}(\rho) = \omega$ and $\alpha$ is $\rho$-stable. Then $\alpha$ is $\omega$-fold stable.
\end{enumerate}
\end{theorem}

\begin{proof}
(1) Let $\alpha_i = f(i)$, where $f$ witnesses $\alpha$'s $\omega$-fold stability, and $\rho = \sup\{\alpha_i: i < \omega\}$. Then by definition we have $L_{\alpha_0} \prec_{\Sigma_1} L_{\alpha_1} \prec_{\Sigma_1} L_{\alpha_2} \prec_{\Sigma_1} \cdots$. By Lemma \ref{StabilityFacts}.2 and \ref{StabilityFacts}.4, we see that $\alpha_i$ is $\alpha_j$-stable whenever $i < j$, and therefore $\alpha_i$ is $\rho$-stable for all $i$. Now let $\tau < \rho$. Then there is some $i < \omega$ so that $\tau < \alpha_i$, and so $\alpha_i$ is $\rho$-stable, witnessing $\rho$'s nonprojectibility in this case.

(2) Let $\{\rho_i: i < \omega\}$ be a cofinal subset of $\rho$, whose existence is guaranteed by our hypothesis $\operatorname{cof}(\rho) = \omega$. For each $i < \omega$, pick a $\delta_i$ so that $\rho_i < \delta_i$ and $\delta_i$ is $\rho$-stable, whose existence is guaranteed by our hypothesis that $\rho$ is nonprojectible. By Lemma \ref{StabilityFacts}.1, $\delta_i$ is $\delta_j$-stable whenever $i < j$. Let $t$ be the least natural number so that $\alpha < \rho_t$. Then define $f: \omega \to \Ord$ by $f(0) = \alpha$ and $f(i+1) = \delta_{i+t}$. Another application of Lemma \ref{StabilityFacts}.1 shows that $\alpha$ is $\delta_i$-stable for all $i \geq t$, and so $f$ witnesses $\alpha$'s $\omega$-fold stability.
\end{proof}

\begin{definition}
\label{Interpretable}
Let $\kappa, \rho$ be ordinals with $\kappa < \rho$. Then $\mathcal{A} \subseteq L_\kappa$ is called $(\kappa, \rho)$-interpretable iff, for any $\Sigma_2$-formula $\varphi$ and parameters $\vec{z}$ with $\mathcal{A} = \{x \in L_\kappa: L_\kappa \models \varphi(x, \vec{z})\}$, if one sets $\mathcal{A}^{L_\rho} = \{x \in L_\rho: L_\rho \models \varphi(x, \vec{z})\}$, then $\mathcal{A} = \mathcal{A}^{L_\rho} \cap L_\kappa$.
\end{definition}

In other words, $\mathcal{A}$ is $(\kappa, \rho)$-interpretable iff any $\Sigma_2$-definition gives a way of extending $\mathcal{A}$ to a subset of $L_\rho$ without adding new elements of $L_\kappa$. Note that the resulting extensions generally depend on the $\varphi$ and $\vec{z}$ chosen, rather than in fact giving unique, canonical extensions.

\begin{lemma}
\label{StabilityToInterpretability}
Let $\kappa, \rho$ be ordinals with $\kappa < \rho$. Then the following are equivalent:

\begin{enumerate}
    \item $\kappa$ is $\Sigma_2$-$\rho$-stable.
    \item Every subset of $L_\kappa$ is $(\kappa, \rho)$-interpretable.
\end{enumerate}
\end{lemma}

\begin{proof}
For the forward direction, assume $\kappa$ is $\Sigma_2$-$\rho$-stable, $\mathcal{A} \subseteq L_\kappa$ and $\varphi, \vec{z}$ are so that $\mathcal{A} = \{x \in L_\kappa: L_\kappa \models \varphi(x, \vec{z})\}$ (if $\mathcal{A}$ is not $\Sigma_2(L_\kappa)$, this is vacuous). By $L_\kappa \prec_{\Sigma_2} L_\rho$, we have $L_\kappa \models \varphi(x, \vec{z})$ iff $L_\rho \models \varphi(x, \vec{z})$ for $x \in L_\kappa$, i.e. $x \in \mathcal{A}$ iff $x \in \mathcal{A}^{L_\rho}$ for all $x \in L_\kappa$. Therefore $\mathcal{A} \cap L_\kappa = \mathcal{A}^{L_\rho} \cap L_\kappa$ and, since $\mathcal{A} \subseteq L_\kappa$, $\mathcal{A} \cap L_\kappa = \mathcal{A}$. Thus, $\mathcal{A} = \mathcal{A}^{L_\rho} \cap L_\kappa$ and, since $\mathcal{A}, \varphi, \vec{z}$ were arbitrary, any subset of $L_\kappa$ is $(\kappa, \rho)$-interpretable.

For the converse direction, let $\vec{z} \in L_\kappa$, and $\varphi$ be a $\Sigma_2$ formula. We aim to show that $L_\kappa \models \varphi(\vec{z})$ iff $L_\rho \models \varphi(\vec{z})$. Let $\mathcal{A} = \{x \in L_\kappa: L_\kappa \models \varphi(\vec{z})\}$. Then, by hypothesis, $\mathcal{A} = \mathcal{A}^{L_\rho} \cap L_\kappa$. Note that, since $\mathcal{A} \subseteq L_\kappa$, we have $\mathcal{A} \cap L_\kappa = \mathcal{A}$. It follows that, for all $x \in L_\kappa$, $x \in \mathcal{A}$ iff $x \in \mathcal{A}^{L_\rho}$, i.e. $L_\kappa \models \varphi(\vec{z})$ iff $L_\rho \models \varphi(\vec{z})$. Since $\varphi(\vec{z})$ is independent of $x$, this gives us the desired result.
\end{proof}

We give our main results now.

\begin{definition}
\label{Sigma2Stability}
Say $\alpha$ is $\Sigma_2$-$\xi$-stable iff $L_\alpha \prec_{\Sigma_2} L_\xi$, where $\prec_{\Sigma_2}$ denotes the relation of being a $\Sigma_2$-elementary substructure.
\end{definition}

So $\Sigma_2$-nonprojectibility is defined analogously, as alluded to previously. The proof of \ref{NonprojectiblePi2Ref} generalizes neatly to show $\Sigma_2$-nonprojectible ordinals are $\Pi_2$-reflecting onto the ordinals $\Sigma_2$-stable up to them.

\begin{theorem}
\label{PredicateElimination}
Let $\kappa < \rho$ be ordinals so that $\kappa$ is $\Sigma_2$-$\rho$-stable and $\mathcal{A} \subseteq L_\rho$ is $\Delta_2^{L_\rho}(L_\kappa)$. Then $\kappa$ is $\mathcal{A}$-$\rho$-stable.
\end{theorem}

\begin{proof}
So $\kappa$ is $\Sigma_2$-$\rho$-stable and $\mathcal{A}$ is $\Delta_2^{L_\rho}(L_\kappa)$, witnessed by a $\Sigma_2$ formula $\varphi$, $\Pi_2$-formula $\psi$ and parameters $\vec{w}, \vec{z} \in L_\kappa$. Then $L_\kappa \models \varphi(x, \vec{w})$ iff $L_\kappa \models \psi(x, \vec{z})$ iff $x \in \mathcal{A}$, for $x \in L_\kappa$, by $\Sigma_2$-stability, so consequently $(L_\kappa, \in, \mathcal{A} \cap L_\kappa) \models \psi(\vec{a})$, where $\psi$ can include the predicate symbol $\mathbf{P}$, iff $(L_\kappa, \in) \models \psi^*(\vec{a})$. Here, $\psi^*$ is obtained from $\psi$ by replacing all positive instances of the predicate symbol $\mathbf{P}(v)$ with $\varphi(v, \vec{z})$ and all negative instances with $\neg \psi(v, \vec{w})$. It is easy to see that this doesn't increase the complexity from $\Sigma_2$ to some higher Lévy rank.

Now, by $\Sigma_2$-stability again, $\mathcal{A}$ is $(\kappa, \rho)$-interpretable and so $(L_\kappa, \in) \models \psi^*(\vec{a})$ iff $(L_\rho, \in) \models \psi^*(\vec{a})$ whenever $\vec{a} \in L_\kappa$. And then the translation works backwards to show that this happens iff $(L_\kappa, \in, \mathcal{A}) \models \psi(\vec{a})$.
\end{proof}

Regarding the formulation of the above theorem, it cannot be improved by much. For example, one can not strengthen this to when $\mathcal{A}$ is $\Delta_2(L_\rho)$. This is because $\{\kappa\}$ is $\Delta_2(L_\rho)$ and $\kappa$ is never $\{\kappa\}$-$\rho$-stable, since $\{\kappa\} \cap L_\kappa = \emptyset$ and $\{\kappa\} \cap L_\rho \neq \emptyset$, thus $L_\kappa$ and $L_\rho$ disagree about the $\Sigma_1$-sentence $\exists x (\mathbf{P}(x))$. This has no bearing on the truth of our main theorem, since we quantify $\mathcal{A}$ before considering the $\mathcal{A}$-$\rho$-stable ordinals. For example, if $\alpha < \kappa$ and $\kappa$ is $\rho$-stable, then $\kappa$ is $\{\alpha\}$-$\rho$-stable. An earlier of this paper had used only $\Sigma_1$-stability, and used $\Sigma_1$ rather than $\Delta_2$. We neglected the possibility that the predicate symbol could occur negatively, raising the complexity from $\Sigma_1$ to $\Sigma_2$ at best. We wish to thank one of the author's friends, who wishes to remain anonymous, for first noticing the mistake, and Philip Welch for discussion regarding whether it could easily be fixed. Part of the argument below is due to him.

For transitive sets $M \subseteq N$, we let $\mathrm{Th}_1^N(M)$ denote the $\Sigma_1$-satisfaction predicate for formulae with parameters in $M$, in $N$; that is, $x \in \mathrm{Th}^N(M)$ iff $x = \langle \ulcorner \varphi \urcorner, \vec{y} \rangle$ for some $\vec{y} \in M$ and $\Sigma_1$-formula $\varphi$ so that $N \models \varphi(\vec{y})$. For example, $M \prec_{\Sigma_1} N$ iff $\mathrm{Th}_1^M(M) = \mathrm{Th}_1^N(M)$.

\begin{lemma}
\label{TruthPredicate}
Let $\tau < \rho$ be ordinals so that $\tau$ is $\mathcal{A}$-$\rho$-stable, where $\mathcal{A} = \mathrm{Th}_1^{L_\rho}(L_\tau)$. Then $\tau$ is $\Sigma_2$-$\rho$-stable.
\end{lemma}

\begin{proof}
By the remark directly before this Lemma, $\mathcal{A} = \mathrm{Th}_1^{L_\tau}(L_\tau)$. Also n.b. that $\mathcal{A} \subseteq L_\tau$, since $L_\tau$ contains all natural numbers and is closed under pairing. Now let $\vec{y} \in L_\tau$ be arbitrary parameters, and $\varphi$ be a $\Sigma_2$-formula. We want to show that $L_\tau \models \varphi(\vec{y})$ iff $L_\rho \models \varphi(\vec{y})$. To achieve this we use the satisfaction predicate. Let $\varphi$ be of the form $\exists x \forall z \psi(z, x, \vec{y})$, for a $\Delta_0$ formula $\psi$. Then:

\begin{equation}
\begin{split}
L_\tau \models \varphi(\vec{y}) & \iff L_\tau \models \exists x \forall z \psi(z, x, \vec{y}) \\ & \iff L_\tau \models \exists x \neg \exists z \neg \psi(z, x, \vec{y}) \\ & \iff \exists x \in L_\tau (\langle \ulcorner \exists z \neg \psi \urcorner, x, \vec{y} \rangle \notin \mathcal{A}) \\ & \iff L_\tau \models \exists x (\neg \mathbf{P}(\langle \ulcorner \exists z \neg \psi \urcorner, x, \vec{y} \rangle)) \\ & \iff L_\rho \models \exists x (\neg \mathbf{P}(\langle \ulcorner \exists z \neg \psi \urcorner, x, \vec{y} \rangle)) \\ & \iff \exists x \in L_\rho (\langle \ulcorner \exists z \neg \psi \urcorner, x, \vec{y} \rangle \notin \mathcal{A}) \\ & \iff L_\rho \models \varphi(\vec{y})
\end{split}
\end{equation}

The idea is that we ``pack'' the $\Pi_1$ part into the $\Delta_0$ part via the truth predicate, then exploit ordinary stability and then unpack the $\Pi_1$ part.
\end{proof}

\begin{corollary}
\label{MainTheorem}
$\rho$ is $\Sigma_2$-nonprojectible iff it is recursively subtle.
\end{corollary}

The forward direction is quite intuitive -- one shows that for sufficiently large $\kappa < \rho$, the additional predicate for $\mathcal{A}$ can be eliminated, by letting $\kappa$ be greater than the ranks of all the parameters in the definition of $\mathcal{A}$. For the converse direction, we use satisfaction predicates, which are sufficiently definable, to get ordinals $\Sigma_2$-stable up to $\rho$.

\begin{proof}
We show the forwards direction first. By Theorem \ref{NonprojectiblePi2Ref} it suffices to show that, for a tail of $\kappa < \rho$, the predicate in the notion of recursive subtlety may be eliminated. Formally restated, we aim to show that, for any $\rho$-recursive $\mathcal{A}$, there is some $\tau < \rho$ so that, for all $\kappa < \rho$ so that $\kappa$ is $\rho$-stable and $\tau < \kappa$, $\kappa$ is $\mathcal{A}$-stable. Then the desired result follows from the fairly obvious fact that if $\alpha$ is $\Pi_n$-reflecting onto $A$ and $\xi < \alpha$, then $\alpha$ is $\Pi_n$-reflecting onto $A \cap [\xi, \alpha)$.

Now let $\mathcal{A}$ be $\rho$-recursively enumerable. By hypothesis, there is some $\Sigma_1$-formula $\varphi$ and parameters $\vec{z} \in L_\rho$ so that $\mathcal{A} = \{x \in L_\rho: L_\rho \models \varphi(x, \vec{z})\}$. For each $0 < i \leq \operatorname{len}(\vec{z})$, let $\alpha_i$ denote the rank in the constructible hierarchy of $z_i$, i.e. $\alpha_i$ is the least ordinal so that $z_i \in L_{\alpha_i + 1}$. Set $\tau = \max\{\alpha_i: 0 < i \leq \operatorname{len}(\vec{z})\}$, and assume that $\kappa < \rho$ is so that $\kappa$ is $\rho$-stable and $\tau < \kappa$. Then, since $\vec{z} \in L_{\tau+1} \subseteq L_\kappa$, it follows that $\mathcal{A}$ is $\Sigma_1^{L_\rho}(L_\kappa)$. By the previous Theorem, $\kappa$ is $\mathcal{A}$-$\rho$-stable. This gives the desired result.

We now show the backwards direction. Suppose $\rho$ is recursively subtle. Let $\kappa < \rho$ be arbitrary. We want to find $\tau < \rho$ so that $\kappa < \tau < \rho$ so that $\tau$ is $\Sigma_2$-$\rho$-stable. Let $\mathcal{A} = \mathrm{Th}_1^{L_\rho}(L_\rho)$. The work of \cite{levy} and \cite{friedman} shows that $\mathcal{A}$ is itself $\rho$-recursively enumerable. By recursive subtlety of $\rho$, there is a $\kappa < \tau < \rho$ so that $\tau$ is $\mathcal{A}$-$\rho$-stable. $\tau$ is as desired, since $\mathrm{Th}_1^{L_\rho}(L_\tau) = \mathcal{A} \cap L_\tau$.
\end{proof}

We would like to note that it is still an open question what recursive analogues of higher large cardinal axioms could be. It is a relatively vague question, and not of such high priority as recursive analogues are typically studied for their applications to recursion theory and proof theory rather than for themselves. However, as a closing remark, we will explain how it is possible that one could consider the recursive analogue of measurability to be $\Sigma_2$-extendibility:

\begin{definition}
\label{Extendible}
An ordinal $\alpha$ is called $\Sigma_2$-extendible iff there is $\beta > \alpha$ so that $\alpha$ is $\Sigma_2$-$\beta$-stable. Let $\zeta$ denote the least $\Sigma_2$-extendible ordinal.
\end{definition}

$\zeta$ does indeed have applications to recursion theory -- for example, the $\Sigma_2(L_\zeta)$ subsets of $\omega$ are precisely the arithmetically quasi-inductive subsets of $\omega$, analogously to how the $\omega_1^{\mathrm{CK}}$-recursive subsets of $\omega$ are precisely the arithmetically inductive subsets of $\omega$. $\zeta$ itself also has a characterisation via generalized computability, namely $\zeta$ is the supremum of the eventually writable ordinals with respect to an infinite time Turing machine. There may be an analogue of hyperarithmetical theory, namely hyperinductive theory, at this stage -- we direct the reader to \cite{klev}.

Say that a sequence $\vec{X} = \langle X_\beta: \beta < \lambda \rangle$ of $\lambda$ many subsets of an admissible ordinal $\kappa > \lambda$ is recursive iff $\{(\alpha, \beta): \alpha \in X_\beta\}$, the subset of $\kappa \times \lambda$ coding $\vec{X}$, is $\kappa$-recursive. It is known that, if $\vec{X}$ is recursive, then $\bigcap_{\beta < \lambda} X_\beta$ is $\kappa$-recursive. Say an ultrafilter $\mathcal{U}$ on $\kappa$ is $L_\kappa$-complete iff, whenever $\vec{X}$ is recursive and $X_\beta \in \mathcal{U}$ for all $\beta < \lambda$, then $\bigcap_{\beta < \lambda} X_\beta \in \mathcal{U}$. Then say $\kappa$ is recursively measurable iff there is an $L_\kappa$-complete nonprincipal ultrafilter on the Boolean algebra of $\kappa$-recursive subsets of $\kappa$.

Then $\kappa$ is recursively measurable iff it is $\Sigma_2$-extendible. This was initially thought to be ``consistent'' in a sense with relations to other large cardinal axioms -- for example, any measurable cardinal is weakly compact and subtle, and:

\begin{proposition}
\label{ExtendibleRefAndNpr}
If $\alpha$ is $\Sigma_2$-extendible, then $\alpha$ is $\Pi_3$-reflecting, nonprojectible, and in fact a limit of smaller $\Pi_3$-reflecting nonprojectible ordinals.
\end{proposition}

The proof is an easy reflection argument, utilizing the truth predicate from \cite{levy}. For more information on the notion of recursive measurability in particular, see \cite{kaufmann}. However, from the work in this paper, we see $\Sigma_2$-extendibility iis likely inadequate, as any measurable cardinal is subtle but $\Sigma_2$-nonprojectible ordinals are limits of many $\Sigma_2$-extendible ordinals. In an upcoming paper, we attempt to tackle recursive analogues in full generality.

\printbibliography[heading=bibintoc,title={References}]

@inbook{rathjen,
    author = "Michael Rathjen",
    title = "The Higher Infinite in Proof Theory",
    series = "Lecture Notes in Logic",
    booktitle = "Logic Colloquium '95: Proceedings of the Annual European Summer Meeting of the Association of Symbolic Logic",
    pages = "275–-304",
    year = "1995",
    DOI = "10.1017/9781316716830.019"
}

@article{rathjen2,
    author = "Michael Rathjen",
    title = "An Ordinal Analysis of Parameter-Free $\Pi^1_2$-Comprehension",
    journal = "Archive for Mathematical Logic",
    volume = "44",
    pages = "263--362",
    year = "2005",
    DOI = "10.1007/s00153-004-0232-4"
}

@misc{jensen,
    author = "Ronald Jensen and Kenneth Kunen",
    title = "Some Combinatorial Properties of L and V",
    howpublished = {Available at \url{https://www.mathematik.hu-berlin.de/~raesch/org/jensen.html}},
    year = "1969"
}

@article{jensen2,
    author = "Ronald Jensen",
    title = "The Fine Structure of the Constructible Hierarchy",
    journal = "Annals of Mathematical Logic",
    volume = "4",
    pages = "229-308",
    year = "1972",
    DOI = "10.1016/0003-4843(72)90001-0"
}

@book{barwise,
    author = "Jon Barwise",
    title = "Admissible Sets and Structures",
    series = "Perspectives in Logic",
    volume = "7",
    year = "1975",
    DOI = "10.1017/9781316717196"
}

@book{simpson,
    author = "Stephen Simpson",
    title = "Subsystems of Second Order Arithmetic",
    series = "Perspectives in Logic",
    year = "2009",
    DOI = "10.1017/CBO9780511581007"
}

@article{arai,
    author = "Toshiyasu Arai",
    title = "Proof Theory for Theories of Ordinals - I: Recursively Mahlo ordinals",
    journal = "Annals of Pure and Applied Logic",
    volume = "122",
    pages = "1--85",
    year = "2003",
    DOI = "10.1016/S0168-0072(03)00020-4"
}

@incollection{richter,
    author = "Wayne Richter and Peter Aczel",
    title = "Inductive Definitions and Reflecting Properties of Admissible Ordinals",
    series = "Studies in Logic and the Foundations of Mathematics",
    booktitle = "Generalized Recursion Theory",
    volume = "79",
    pages = "301--381",
    year = "1974",
    DOI = "10.1016/S0049-237X(08)70592-5"
}

@book{levy,
    author = "Azriel Lévy",
    title = "A Hierarchy of Formulas in Set Theory",
    series = "Memoirs of the American Mathematical Society",
    year = "1965",
    DOI = "10.2307/2270349"
}

@mastersthesis{klev,
    author = "Ansten Mørch-Klev",
    title = "Extending Kleene’s $\mathcal{O}$ Using Infinite Time Turing Machines, or How With Time She Grew Taller and Fatter",
    school = "Institute of Logic, Language and Computation, Universiteit van Amsterdam",
    year = "2007"
}

@inbook{kaufmann,
    author = "Matt Kaufmann",
    title = "On Existence of $\Sigma_n$ End Extensions",
    series = "Lecture Notes in Mathematics",
    booktitle = "Logic Year 1979--80",
    pages = "92--103",
    year = "1981",
    DOI = "10.1007/BFb0090942"
}

@article{friedman,
    author = "Sy-David Friedman",
    title = "The $\Sigma^*$ Approach to the Fine Structure of $L$",
    journal = "Fundamenta Mathematicae",
    volume = "154",
    pages = "133--158",
    year = "1997",
}
\end{document}